\documentclass[12pt]{article}
\usepackage{amsgen,amsmath,amstext,amsbsy,amsopn,amsfonts,amssymb}
\newtheorem{theorem}{Theorem}

\newtheorem{proposition}[theorem]{Proposition}
\newtheorem{obs}[theorem]{Observation}

 \newtheorem{defi}[theorem]{Definition}

\newtheorem{exa}[theorem]{Example}

\newtheorem{rem}[theorem]{Remark}

\newtheorem{rems}[theorem]{Remarks}

\newtheorem{ack}[theorem]{Acknowlegment}

\def\bsq{\blacksquare\medskip}

\def\H{\mathcal H}

\def\CCC{{\mathbb C}}

\def\RR+{{\mathbb R}^*}

\def\Q_p{{\mathbb Q}_p}

\def\eps{\varepsilon}
\def\Ga{\Gamma}

\begin{document}

\title{Lattices with  and lattices without spectral gap}
\author{Bachir Bekka and Alexander Lubotzky
 \footnote{This research was supported by a grant from the ERC}}
\date{September 10, 2009 }

\maketitle

\rightline{\it For  Fritz Grunewald on his 60th birthday}

\begin{abstract}
Let $G={\mathbf G}(\mathbf k)$ be the $\mathbf k$-rational points
of a simple algebraic group ${\mathbf G}$ over  a local field $\mathbf k$
and let $\Gamma$ be  a lattice  in $G.$ We show that the  regular representation  $\rho_{\Ga\backslash G}$
of $G$  on $L^2(\Ga\backslash G)$ has a spectral gap, that is, the restriction of $\rho_{\Ga\backslash G}$ to the
orthogonal of the constants in $L^2(\Ga\backslash G)$
has no almost  invariant vectors.
On the other hand, we give examples of  locally compact simple groups $G$ and lattices $\Ga$
for which  $L^2(\Ga\backslash G)$ has  no spectral gap. 
This answers  in the negative    a question asked by Margulis \cite[Chapter III, 1.12]{Margulis}.
In fact, $G$ can be taken to be the 
 group of orientation preserving automorphisms of a $k$-regular tree for $k>2.$

\end{abstract}

\section{Introduction}
\label{S1}
Let $G$ be a locally compact  group. 
Recall that a unitary representation $\pi$
of $G$  on a Hilbert space $\cal H$   has almost invariant vectors 
if, for every compact subset $Q$ of $G$ and every 
$\varepsilon>0,$ there exists a unit vector $\xi\in\H$ such that
$\sup_{x\in Q}\Vert \pi(x)\xi-\xi\Vert <\varepsilon.$
If this holds, we also say that the trivial representation $1_G$
is weakly contained in $\pi.$ 

Recall that a lattice $\Ga$ in $G$ is a discrete subgroup
such that  there exists  a finite $G$-invariant regular Borel measure
$\mu$ on $\Ga\backslash G$. 
Denote by $\rho_{\Ga\backslash G}$ the unitary representation
of $G$ given by 
right translation on the Hilbert space $L^2(\Ga\backslash G,\mu)$ 
of the square integrable measurable functions
on $\Ga\backslash G.$  The  subspace $\CCC 1_{\Ga\backslash G}$ of the constant functions on $\Ga\backslash G$ is 
$G$-invariant as well as its orthogonal
complement 
$$
L^2_0(\Ga\backslash G)=\left\{\xi\in L^2(\Ga\backslash G)\, :\, \int_{\Ga\backslash G} \xi(x) d\mu(x)=0\right\}.
$$
Denote by $\rho_{\Ga\backslash G}^0$ the restriction of $\rho_{\Ga\backslash G}$
to $L^2_0(\Ga\backslash G,\mu)$.
We say that $\rho_{\Ga\backslash G}$ (or $L^2(\Ga\backslash G,\mu)$)  has a \emph{spectral gap}  if  $\rho_{\Ga\backslash G}^0$   has  no  almost invariant vectors.
(In  \cite[Chapter III., 1.8]{Margulis}, $\Ga$ is then called weakly cocompact.)
It is well-known  that  $L^2(\Ga\backslash G)$ has a spectral gap
when $\Ga$ is cocompact in $G$ 
(see \cite[Chapter III,  1.10]{Margulis}).
Margulis  ({\it op.cit}, 1.12) asks whether this result holds  
more generally when $\Gamma$ 
is a subgroup of finite covolume.

The goal of this note is to prove the following results:

\begin{theorem}
\label{Theo1} Let $\mathbf G$ be a simple algebraic  group over  a local field $\mathbf k$
and $G=\mathbf G(\mathbf k),$ the group of $\mathbf k$-rational points in $\mathbf G$ . Let
$\Ga$ be a  lattice in $G$.
Then the  unitary representation $\rho_{\Ga\backslash G}$ on $L^2(\Ga\backslash G)$ has a spectral gap.
\end{theorem}
\begin{theorem}
\label{Theo2} For an integer $k>2,$ let $X$ be the $k$--regular tree and  $G={\rm Aut} (X).$
Then $G$ contains a lattice $\Ga$ for which  the unitary representation $\rho_{\Ga\backslash G}$ on $L^2(\Ga\backslash G)$ has no spectral gap.
\end{theorem}

So, Theorem \ref{Theo2} answers in the negative Margulis' question mentioned above.

Theorem~\ref{Theo1} is known in case $\mathbf k=\mathbf R$ 
(\cite{Bekka}).  It  holds, more generally,  when $G$ is a real Lie group (\cite{Bekka-Cornulier}).
Observe also that when  $\mathbf k-{\rm rank} (\mathbf G)\geq 2,$ the group 
$G$ has Kazhdan's Property (T) (see \cite{BHV}) and Theorem \ref{Theo1} is clear in this case.
When $\mathbf k$ is non-archimedean 
with characteristic $0$, every lattice $\Ga$ in $\mathbf G(\mathbf k)$
is uniform (see \cite[p.84]{Serre}) and hence the result
holds as mentioned above. By way of contrast, $G$ has many 
non uniform lattices when the characteristic of $\mathbf k$ is non zero
(see \cite{Serre} and \cite{Alex-rank1}).
So, in order to prove Theorem~\ref{Theo1}, it suffices to consider the case where the characteristic of $\mathbf k$ is non-zero
and where $\mathbf k-{\rm rank} (\mathbf G)=1.$ 

Recall that when $\mathbf k$ is non-archimedean  and $\mathbf k-{\rm rank} (\mathbf G)=1,$ the group $\mathbf G(\mathbf k)$ acts by automorphisms on the associated Bruhat-Tits tree $X$ (see \cite{Serre}).
This tree is either the $k$-regular tree $X_k$ (in which every vertex has constant degree $k$)
or is the bi-partite bi-regular tree $X_{k_0,k_1}$ (where every vertex has either degree $k_0$ or degree $k_1$
and where all neighbours of a vertex of degree $k_{i}$ have degree $k_{1-i}$  ).
The proof of Theorem~\ref{Theo1} will use the special structure  of a fundamental
domain for the action of $\Ga$ on $X$ as described in \cite{Alex-rank1} (see
also \cite{Raghunathan-FundamentalDomain} and \cite{Baum}).

 Theorems  \ref{Theo1} and \ref{Theo2}  provide a further illustration
of   the  different behaviour of  general tree lattices 
as compared to  lattices in rank one simple Lie groups over local fields;
for more on this topic, see \cite{AlexEdin}.

%However, as Theorem~\ref{Theo2} shows, it does not hold for an arbitrary locally compact group $G.$

The proofs of Theorems~\ref{Theo1} and \ref{Theo2} will be given in Sections~\ref{S3} and \ref{S4};
they rely in a crucial way on Proposition~\ref{CorExpander} 
 from Section~\ref{S2},  which relates the existence of a spectral gap with
expander diagrams. In turn,  Proposition~\ref{CorExpander}
is based, much in the spirit of \cite{Brooks},
on analogues for diagrams proved in \cite{Mokhtari} and \cite{Morgenstern1}
 of the inequalities of   Cheeger and  Buser between
the isoperimeric constant  and the bottom 
of the spectrum of the Laplace operator on a Riemannian manifold
(see Proposition~\ref{TheoExpander}).  This connection between the combinatorial expanding 
property and representation theory is by now a very popular theme;  see
\cite{Alex-book} and the references therein. While most applications in this monograph
are from representation theory to combinatorics, we use in the current paper
this connection in the opposite direction:  the existence or absence of a spectral gap
is deduced from the existence of an expanding diagram or of  a non-expanding diagram,
respectively.

\section{Spectral gap and expander diagrams}
\label{S2}
We first show how the existence of a spectral gap
for groups acting on trees is related with the 
bottom of the spectrum
of the Laplacian for an associated diagram.

A graph $X$ consists of  a set of vertices $VX$, a set of oriented edges $EX$,
a fix-point free involution $^-:EX\to EX,$ and end point mappings $\partial_i: EX\to VX$
for $i=0,1$ such that $\partial_i(\overline e)= \partial_{1-i}(e)$ for all $e\in EX.$
Assume that $X$ is  locally finite, that is, for every $x\in VX$, the degree ${\rm \deg} (x)$ of 
$x$ is finite, where ${\rm \deg }(x)$ is the cardinality of the set 
$$\partial_0^{-1}(x)=\{e\in EX\ : \ \partial_0(e)=x\}.$$
  The group ${\rm Aut} (X)$ of automorphisms
of the graph $X$ is a locally compact group in the topology of pointwise convergence
on $X,$ for which the stabilizers of vertices are compact  open subgroups.

We will   consider infinite graphs  called diagrams of finite volume. 
An \emph{edge-indexed graph} $(D,i)$ is a graph $D$ equipped 
with a function $i:ED\to \mathbf R^+$ (see \cite[Chapter 2]{BL}).
A measure $\mu$ for an edge-indexed graph $(D,i)$
is a function $\mu: VD\cup ED\to \mathbf R^+$
with the following properties (see \cite{Mokhtari} and \cite[2.6]{BL}):
\begin{itemize}
\item $i(e)\mu(\partial_0 e)=\mu(e)$
\item $\mu(e)=\mu(\overline e)$  for all $e\in VD$, and
\item $\sum_{x\in VD} \mu(x)<\infty.$
\end{itemize}
Following \cite{Morgenstern1}, we will say that $D=(D,i,\mu)$ is a \emph{diagram of finite volume.}
The in-degree ${\rm indeg}(x)$ of  a vertex $x\in VD$ is defined by
$${\rm indeg}(x)= \sum_{e\in\partial_0^{-1}(x)} i(e)=\sum_{e\in\partial_0^{-1}(x)} \frac{\mu(e)}{\mu(x)}.$$
The diagram $D$ is $k$-regular if ${\rm indeg}(x)=k$ for all $x\in VD.$

Let $D=(D,i,\mu)$ be a connected diagram of finite volume. Observe that 
$\mu$ is determined, up to a multiplicative constant, by the weight function $i.$
Indeed, fix $x_0 \in VD$ and set $\Delta(e)= i(e)/i(\overline e)$ for $e\in ED.$
Then 
$$
\mu(\partial_1 e)= \frac{\mu(\overline e)}{i({\overline e})}= \frac{\mu(e)}{i(\overline e)}=\mu(\partial_0 e)\Delta(e)
$$ for every $e\in ED.$ Hence
$\mu(x)=\Delta(e_1)\Delta(e_2)\dots \Delta(e_n)\mu(x_0)$ 
for every path $(e_1,e_2,\dots, e_n)$ from $x_0$ to $x\in VD.$

Let $D=(D,i,\mu)$ be a diagram of finite volume. An inner product is defined for functions 
on $VD$ by 
$$\langle f, g \rangle= \sum_{x\in VD} f(x)\overline{g(x)} \mu(x).$$
The Laplace operator $\Delta$ on functions $f$  on $VD$  is defined by 
$$
\Delta f(x) = f(x)- \frac{1}{{\rm indeg}(x)}\sum_{e\in\partial_0^{-1}(x)} \frac{\mu(e)}{\mu(x)}f(\partial_1 (e)).
$$
%where the in-degree of  a vertex $x$ in $D$ is 
%$${\rm indeg}(x)= \sum_{(x,y)\in ED} i(x,y)=\sum_{(x,y)\in ED} \frac{\mu(x,y)}{\mu(x)}.$$
The operator $\Delta$ is a self-adjoint positive operator on 
$L^2(VD).$
Let 
$$L^2_0(VD)=\{f\in L^2(VD)\ :\ \langle f, 1_{VD} \rangle=0\}$$
and set
$$
\lambda(D)= \inf_{f} \langle \Delta f, f \rangle,$$
where $f$ runs over the  unit sphere in $L^2_0(VD)$.
Observe that 
$$
\lambda(D)= \inf \{ \lambda\ :\ \lambda\in\sigma (\Delta)\setminus\{0\}\},
$$
where $\sigma (\Delta)$ is the spectrum of $\Delta.$

Let now $X$ be a locally finite tree, and let $G$ be a closed subgroup of ${\rm Aut} (X)$.
Assume that $G$ acts with finitely many orbits on $X$.  
Let $\Ga$ be a discrete subgroup of $G$ acting without inversion on $X.$
Then the quotient graph  $\Ga\backslash X$ is well-defined.
Since $\Ga$ is discrete, for every vertex $x$ and every edge $e,$ the  stabilizers $\Ga_x$ and $\Ga_e$ 
are finite. Moreover, $\Ga$ is a lattice in $G$ if and only if
$\Ga$ is a lattice in ${\rm Aut} (X)$ and this happens if and only if 
$$\sum_{x\in D} \frac{1}{|\Ga_x| }<\infty,$$
where $D$ is a fundamental domain of $\Ga$ in $X$ (see \cite{Serre}).
The quotient graph  $\Ga\backslash X\cong D$ is endowed with the structure of
an  edge-indexed graph  given by the weight function $i: ED\to \mathbf R^+$ where
$i(e)$ is the index of $\Ga_e$ in $\Ga_x$ for $x=\partial_0(e)$.
A measure 
$\mu: VD\cup ED\to \mathbf R^+$ is defined by 
$$
\mu(x)=  \frac{1}{|\Ga_x|} \qquad \text{and} \qquad \mu(e)=  \frac{1}{|\Ga_e|}
$$
for $x\in VD$ and $e\in ED$. Observe that $\mu(VD)=\sum_{x\in D} {1}/{|\Ga_x|} <\infty.$
So, $D=(D,i,\mu)$ is a diagram of finite volume.
\medskip

Let $G$ be a group acting on a tree $X.$ As in \cite[0.2]{BurgerMozes}, we say that the action of $G$
on $X$ is \emph{locally $\infty$-transitive} if, for every $x\in VX$ and every $n\geq 1,$ the stabilizer $G_x$ of $x$  acts transitively  on the sphere $\{y\in X\ : d(x,y)= n\}.$

\begin{proposition} 
\label{proposition1} Let $X$ be either the $k$-regular  tree $X_k$ 
or the  bi-partite bi-regular tree $X_{k_0,k_1}$
for $k\geq 3$ or $k_0\geq3$ and $k_1\geq 3.$
Let $G$ be a closed subgroup of ${\rm Aut} (X).$ Assume that
the following conditions are both satisfied:
\begin{itemize}
 \item $G$ acts transitively on $VX$ in the case $X=X_k$ and 
$G$ acts transitively on the set of vertices of degree $k_0$ as well as on 
  the set of vertices of degree $k_1$ in the case $X=X_{k_0,k_1};$
\item the action of $G$ on $X$ is locally $\infty$-transitive.
\end{itemize}
Let $\Ga$ be a lattice in $G$ and let $D= \Ga\backslash X$ be the
corresponding diagram of finite volume.
 The following properties are equivalent:
\begin{itemize}
 \item[(i)] the  unitary representation $\rho_{\Ga\backslash G}$ on $L^2(\Ga\backslash G)$ has a spectral gap;
\item[(ii)] $\lambda(D)>0.$
\end{itemize}
 \end{proposition}

 For the proof of this proposition, we will need
a few  general facts.
Let $G$ be a second countable  locally compact group  and $U$ a compact  subgroup of $G.$
Let $C_c(U\backslash G/U)$ be the space  of  continuous functions $f:G\to \mathbf C$
which have compact support and which are constant on the double cosets $UgU$ for $g\in G.$ 

Fix a left Haar measure $\mu$ on $G.$ 
Recall  that   $L^1(G,\mu) $ is a Banach algebra 
under the convolution product, the $L^1$-norm
and the involution $f^*(g)= \overline{f(g^{-1})}$;
observe  that   $C_c(U\backslash G/U)$
is a  $*$-subalgebra of $L^1(G,\mu).$
Let $\pi$ be a (strongly continuous) unitary representation
of $G$ on a Hilbert space $\cal H.$ 
A continuous $*$-representation of $L^1(G)$, still denoted by $\pi,$
is defined on $\cal H$ by 
$$\pi(f)\xi= \int_G f(x) \pi(x)\xi d\mu(x),\qquad f\in L^1(G),\quad \xi\in\cal H.$$
Assume that the closed subspace
${\cal H}^U$ of $U$-invariant vectors in $\H$ is non-zero.
Then $\pi(f){\cal H}^U\subset {\cal H}^U$ 
for all $f\in C_c(U\backslash G/U)$.
In this way,   a  continuous $*$-representation $\pi_U$
 of $C_c(U\backslash G/U)$ is defined on ${\cal H}^U$.

\begin{proposition}
\label{proposition3}
 With the previous notation,
let $f\in C_c(U\backslash G/U)$ be a function 
with the following properties: $f(x)\geq 0$ for all $x\in G,$   $\int_Gf d\mu =1,$ and 
 the subgroup generated by the support of $ f$ is dense in $G.$
The following conditions are equivalent:
\begin{itemize}
 \item [(i)] the trivial representation $1_G$
is  weakly contained in $\pi;$ 
\item [(ii)] $1$ belongs to the spectrum of the operator  $\pi_U(f).$ 
\end{itemize}
\end{proposition}
\begin{proof}
Assume that $1_G$
is   weakly contained in $\pi.$  There exists a sequence of unit vectors
$\xi_n\in\H$ such that
$$
\lim_n\Vert \pi(x)\xi_n-\xi_n\Vert=0,
$$
 uniformly  over compact subsets of $G.$ Let 
$$\eta_n=\int_U \pi(u) \xi_n du,$$
where $du$ denotes the normalized Haar measure on $U.$
It is easily checked that  $\eta_n\in {\cal H}^U$ and that
$$\lim_n\Vert \pi(f)\eta_n-\eta_n\Vert=0.$$
Since
$$
\Vert\eta_n -\xi_n\Vert \leq \int_U \Vert \pi(u) \xi_n -\xi_n\Vert du,
$$
we have $\Vert \eta_n\Vert \geq 1/2$ for sufficiently large $n.$
This shows that $1$ belongs to the spectrum of the  operator 
$\pi_U(f).$ 

For the converse, assume that $1$ belongs to the spectrum of 
$\pi_U(f).$ Hence, $1$ belongs to the spectrum of 
$\pi(f),$ since $\pi_U(f)$ is the restriction of 
$\pi(f)$ to the invariant subspace ${\cal H}^U.$
As the subgroup generated by the support of $f$ is dense in $G$, this   implies
that $1_G$ is weakly contained in $\pi$ (see \cite[Proposition G.4.2]{BHV}).
\end{proof}
\bigskip

\noindent
\textbf{Proof of Proposition~\ref{proposition1}  }
We give the proof  only in the case where  $X$ is the  bi-regular tree $X_{k_0,k_1}.$
The case where $X$ is  the regular tree $X_k$ is similar and even simpler.

Let $X_0$ and $X_1$ be the subsets of $X$ consisting of the vertices
of degree $k_0$ and $k_1,$ respectively.
Fix two  points $x_0\in X_0$ and $x_1\in X_1$ with $d(x_0, x_1)=1.$
So, $X_0$ is the set of vertices  $x$ for which  $d(x_0,x)$ is even and
$X_1$ is the set of vertices  $x$ for which $d(x_0,x)$ is odd.
Let $U_0$ and $U_1$ be the stabilizers of $x_0$ and $x_1$ in $G.$
 Since $G$ acts transitively on $X_0$ and on $X_1,$ we have $G/U_0\cong X_0$.
and $G/U_1\cong X_1$.

We can view  the normed  $*$-algebra $C_c(U_0\backslash G/U_0)$
  as a space of finitely supported functions on $X_0$.
Since $U_0$ acts transitively on every sphere around $x_0$, 
it is well-known that the pair $(G,U_0)$ is a Gelfand pair,
that is, the algebra $C_c(U_0\backslash G/U_0)$ is commutative
(see for instance \cite[Lemma 2.1]{BLRW}).
Observe that $C_c(U_0\backslash G/U_0)$ is the linear span of  the
characteristic functions $\delta_n^{(0)}$ (lifted to $G$) of  spheres
of even radius  $n$ around $x_0$.  
Moreover, $C_c(U_0\backslash G/U_0)$ is generated by $\delta_2^{(0)}$; indeed,
 this follows from the formulas (see  \cite[Theorem 3.3]{BLRW})
\begin{align*}
\delta_{4}^{(0)}&= \delta_{2}^{(0)}\ast\delta_{2}^{(0)}-k_0(k_1-1)\delta_0^{(0)} -(k_1-2)\delta_{2}^{(0)}\\
\delta_{2n+2}^{(0)}&= \delta_{2}^{(0)}\ast\delta_{2n}^{(0)}-(k_0-1) (k_1-1)\delta_{2n-2}^{(0)} -(k_1-2)\delta_{2n}^{(0)} \qquad \text{for}\quad n\geq 2.
\end{align*}
Let $f_0=\dfrac{1}{\Vert\delta_2^{(0)}\Vert_1} \delta_2^{(0)}$.
We claim that $f_0$ has all the properties listed in Proposition~\ref{proposition3}.

Indeed, $f_0$ is a non-negative and $U_0$-bi-invariant function on $G$ with  $\int_G f_0(x) dx=1.$
Moreover,   let $H$ be the closure of  the subgroup generated by the support of $f_0$.
Assume, by contradiction, that $H\neq G.$ Then there exists a function  in $ C_c(U_0\backslash G/U_0)$
whose support is disjoint from $H.$  This is a contradiction, as the algebra $C_c(U_0\backslash G/U_0)$ is generated by 
$f_0.$ This shows that $H=G$.

Let $\pi$ be  the unitary representation of $G$ on $L_0^2(\Ga\backslash G)$
defined by right translations.  Observe that the space of $\pi(U_0)$-invariant  vectors
is $L_0^2(\Ga\backslash X_0)$.
So, we have a $*$-representation $\pi_{U_0}$ of $C_c(U_0\backslash G/U_0)$
on  $L^2(\Ga\backslash X_0, \mu)$, where $\mu$ is the measure
on the diagram $D= \Ga\backslash X,$ as defined above.
%Let $\pi_{U_0}$ denote the restriction of this representation to the  $C_c(U_0\backslash G/U_0)$-invariant subspace  $$L^2_0(\Ga\backslash X_0, \mu)
%=\{\xi \in L^2(\Ga\backslash X_0, \mu)\ :\ \sum_{x\in X_0} \xi(x)\mu(x)=0 \}$$

Similar facts are also true for the algebra $C_c(U_1\backslash G/U_1):$
this is a commutative normed $*$-algebra, it is generated by the characteristic function $\delta_2^{(1)}$ of the sphere of radius 2 around $x_1$, and
the representation $\pi$ of $G$
on $L_0^2(\Ga\backslash G)$ induces a   $*$-representation $\pi_{U_1}$  of $C_c(U_1\backslash G/U_1)$
on $L_0^2(\Ga\backslash X_1, \mu).$ 
Likewise, the function $f_1=\dfrac{1}{\Vert\delta_2^{(1)}\Vert_1} \delta_2^{(1)}$
 has all the properties listed in Proposition~\ref{proposition3}.

%Let $\Delta$ be the Laplace operator on $D.$
%We claim  that $\lambda(\Delta) >0$ if and only if $1$ does not belong to the spectra of the operators $\pi_{U_0} (f_0)$ and $\pi_{U_1}(f_1)$.
%Once proved, Proposition~\ref{proposition1}  follows from Proposition~\ref{proposition3}.
Let $A_X$ be the adjacency operator defined on $\ell^2(X)$
by 
$$
A_Xf(x)= 
\frac{1}{{\rm deg}(x)}\sum_{e\in\partial_0^{-1}(x)} f(\partial_1(e)),
\qquad f\in \ell^2(X).
$$
 Since $A_X$ commutes with
automorphisms of $X$, it
induces an operator $A_D$  on $L^2(VD, \mu)$
given by 
$$
A_Df(x)= \frac{1}{{\rm indeg}(x)}\sum_{e\in\partial_0^{-1}(x)} \frac{\mu(e)}{\mu(x)}f(\partial_1(e)),
\qquad f\in L^2(VD,\mu), 
$$
where $D$ is the  diagram obtained from the quotient graph $\Ga\backslash X.$
So, $\Delta= I-A_D,$ where $\Delta$ is the Laplace operator on $D.$

Let  $B_D$ denote the restriction of $A_D $ to the space $L^2_0(VD,\mu).$ 
It follows that $\lambda(\Delta) >0$ if and only if $1$ does not belong to the spectrum 
of $B_D.$

 Proposition~\ref{proposition1} will be proved,  once we have shown the following
 
 \medskip
\noindent
{\bf Claim:} $1$  belongs to the spectrum 
of  $B_D $  if and only if 
$1_G$ is weakly contained in $\pi.$

\medskip

For this, we  consider the squares of the operators $A_X$ and $A_D$ and compute
$$
A_X^2f(x)= 
\frac{1}{ k_0k_1}{\rm deg}(x) f(x) + \frac{1}{ k_0k_1} \sum_{d(x,y)=2} f(y),\qquad f\in \ell^2(X).
$$
The subspaces $\ell^2(X_0)$ and $\ell^2(X_1)$ of $\ell^2(X)$ are invariant under $A_X^2$
and the restrictions of $A_X^2$ to  $\ell^2(X_0)$ and $\ell^2(X_1)$ are given
by right convolution with  the functions 
\begin{align*}
 g_0&= \frac{1}{ k_0k_1}\delta_{e} + (1- \frac{1}{ k_0k_1}) f_0 \\
g_1&= \frac{1}{ k_0k_1}\delta_{e} + (1- \frac{1}{ k_0k_1}) f_1, \\
\end{align*}
where $\delta_e$ is the Dirac function at the group unit $e$ of $G.$

It follows that the restrictions of $B_D^2$ to the subspaces $L^2_0(\Ga\backslash X_0, \mu)$
and $L^2_0(\Ga\backslash X_1, \mu)$ coincide with the operators 
$\pi_{U_0}(g_0)$ and $\pi_{U_1}(g_1),$ respectively.

For $i=0,1,$ the spectrum $\sigma (\pi_{U_i}(g_i))$ 
 of $\pi_{U_i}(g_i)$ 
is the set 
$$
\sigma (\pi_{U_i}(g_i))=\left\{\frac{1}{ k_0k_1} + (1- \frac{1}{ k_0k_1}) \lambda\ : \lambda\in \sigma (\pi_{U_i}(f_i)) \right\}.
$$
Thus, $1$  belongs to the spectrum of  $\pi_{U_0} (f_i)$ 
if and only if $1$  belongs to the spectrum of  $\pi_{U_0} (g_i)$.

To prove the claim above, assume that  $1$ belongs to the spectrum of $B_D.$ Then $1$ belongs to the spectrum of $B_D^2.$
Hence $1$ belongs to the spectrum of either $\pi_{U_0}(g_0)$ or $\pi_{U_1}(g_1)$
and therefore $1$ belongs to the spectrum of either $\pi_{U_0}(f_0)$ or $\pi_{U_1}(f_1).$
It follows from Proposition~\ref{proposition3} that $1_G$ is weakly contained in $\pi.$

Conversely,  suppose that $1_G$ is weakly contained in $\pi$.
Then, again by Proposition~\ref{proposition3}, $1$ belongs to the spectra of  $\pi_{U_0}(f_0)$ and $\pi_{U_1}(f_1).$ Hence,  $1$ belongs to the spectra of  $\pi_{U_0}(g_0)$ and $\pi_{U_1}(g_1).$
We claim that $1$  belongs to the spectrum of $B_D.$ 

Indeed, 
assume by contradiction that $1$  does not belong to the spectrum of $B_D,$
that is, $B_D-I$ has a bounded inverse on $L^2_0(VD, \mu).$
Since $1$ belongs to the spectrum of  the self-adjoint operator $\pi_{U_0}(g_0),$
there exists  a sequence   of unit vectors $\xi_n^{(0)}$
in $L^2_0(\Ga\backslash X_0, \mu)$ with
$$
 \lim_n\Vert\pi_{U_0}(g_0)\xi_n^{(0)}-\xi_n^{(0)}\Vert =0. 
%\qquad\text{and}\qquad \lim_n\Vert\pi_{U_1}(g_1)\xi_n^{(1)}-\xi_n^{(0)}\Vert =0.
$$
As the restriction of  $B_D^2$ to  $L^2_0(\Ga\backslash X_0, \mu)$
 coincides with 
$\pi_{U_0}(g_0)$, we have 
\begin{align*}
\Vert\pi_{U_0}(g_0)\xi_n^{(0)}-\xi_n^{(0)}\Vert &=\Vert (B_D^2-I)\xi_n^{(0)}\Vert\\
&=\Vert (B_D-I)(B_D+I)\xi_n^{(0)}\Vert\\
&\geq \frac{1}{\Vert (B_D -I)^{-1}\Vert}\Vert(B_D+I)\xi_n^{(0)}\Vert
\end{align*}
So, $\lim_n \Vert B_D\xi_n^{(0)} +\xi_n^{(0)}\Vert=0$.
On the other hand, observe that $B_D$ maps $L^2_0(\Ga\backslash X_0, \mu)$ to the subspace
$L^2(\Ga\backslash X_1, \mu)$ and that these subspaces are orthogonal to each other.
Hence, 
$$
\Vert B_D\xi_n^{(0)} +\xi_n^{(0)}\Vert^2=\Vert B_D\xi_n^{(0)}\Vert^2 +\Vert\xi_n^{(0)}\Vert^2
$$ 
This is a contradiction since $\Vert \xi_n^{(0)}\Vert=1$ for all $n.$
The proof of Proposition~\ref{proposition1} is now complete.$\bsq$

\bigskip
Next, we rephrase Proposition~\ref{proposition1} in terms of expander diagrams.
Let $(D, i,w)$ be a  diagram with finite volume.
For a subset $S$ of $VD,$ set 
$$E(S, S^c)= \{e\in ED\ :\ \partial_0(e)\in S,\  \partial_1(e)\notin S \}.$$
We say that $D$ is an \emph{expander diagram} if there exists $\eps>0$ such that
$$
\frac{\mu (E(S, S^c))}{\mu(S)}\geq \eps
$$
for all $S\subset VD$ with $\mu(S)\leq \mu (D)/2.$
The motivation for this definition
comes from expander graphs (see \cite{Alex-book}). 

We quote from  \cite{Mokhtari} and \cite{Morgenstern1}  the following result  
which is standard in the case of finite graphs.

\begin{proposition} \textbf{(\cite{Mokhtari}, \cite{Morgenstern1})}
\label{TheoExpander}
Let $(D,i, w)$ be a  diagram with finite volume.
Assume that $\sup_{e\in ED} i(\overline e)/i(e) <\infty$ and that
$\sup_{x\in VD} {\rm indeg} (x)<\infty$
  The following
conditions are equivalent:
\begin{itemize}
 \item [(i)]$D$ is an expander diagram; 
\item [(ii)] $\lambda(D)>0.$
\end{itemize}

\end{proposition}
As an immediate consequence of Propositions~\ref{proposition1}
and \ref{TheoExpander},
we obtain the following result which relates
the existence of a spectral gap to an expanding property of
the corresponding diagram. 

\begin{proposition}
 \label{CorExpander}
Let $X$ be either the $k$-regular  tree $X_k$ 
or the  bi-partite bi-regular tree $X_{k_0,k_1}$
for $k\geq 3$ or $k_0\geq3$ and $k_1\geq 3.$
Let $G$ be a closed subgroup of ${\rm Aut} (X)$
satisfiying both conditions from Proposition~\ref{proposition1}.
Let $\Ga$ be a lattice in $G$ and let $D= \Ga\backslash X$ be the
corresponding diagram of finite volume.
 The following properties are equivalent.
\begin{itemize}
 \item[(i)] The  unitary representation $\rho_{\Ga\backslash G}$ on $L^2(\Ga\backslash G)$ has a spectral gap;
\item[(ii)] $D$ is an expander diagram.
\end{itemize}
\end{proposition}

.

\section{Proof of Theorem~\ref{Theo1}}
\label{S3}
Let $G={\mathbf G}(\mathbf k)$ be the $\mathbf k$-rational points
of a simple algebraic group ${\mathbf G}$ over  a local field $\mathbf k$
and let $\Gamma$ be  a lattice  in $G.$ As explained
in the Introduction, we may assume that 
$\mathbf k$ is non-archimedean and that $\mathbf k-{\rm rank} (\mathbf G)=1.$
By the Bruhat-Tits theory, $G$ acts on a regular or bi-partite bi-regular tree $X$ with one or two orbits.
Moreover, the action of $G$ on $X$ is locally $\infty$-transitive
(see \cite[p.33]{Choucroun}).

Passing to the subgroup $G^+$ 
of index at most two consisting of orientation preserving automorphisms, we can assume that $G$ acts
without inversion. Indeed, assume that $L^2(\Ga\cap G^+\backslash G^+)$ has a spectral gap.
If $\Ga$ is contained in $G^+,$ then $L^2(\Ga\backslash G)$ has a spectral gap since $G^+$ has
finite index (see \cite[Proposition 6]{Bekka-Cornulier}).
If $\Ga$ is not contained in $G^+,$ then $\Ga\cap G^+\backslash G^+$ may be identified as a $G^+$-space with $\Ga\backslash \Ga G^+=\Ga\backslash G.$ Hence, $1_{G^+}$ is not weakly contained in the $G^+$-representation 
defined on  $L^2_0(\Ga\backslash G).$

Let $X$ be the Bruhat-Tits tree associated to $G.$
 It is shown in
\cite[Theorem 6.1]{Alex-rank1} (see also
\cite{Baum}) that $\Ga$ has  fundamental domain  $D$ in $X$ of the following form:
there exists a finite set $F\subset D$ such that $D\setminus F$ is a union
of finitely many disjoint rays $r_1,\dots, r_s$. (Recall that a ray in $X$ is an infinite 
path beginning at some vertex and without backtracking.)
Moreover, for every ray $r_j=\{x_0^j,x_1^j,x_2^j, \dots \}$ in $D\setminus F$,
the stabilizer $\Ga_{x_{i}^j}$ of $x_{i}^j$ is contained
in the stabilizer $\Ga_{x_{i+1}^j}$ of $x_{i+1}^j$ for all $i.$

To prove  Theorem~\ref{Theo1}, we apply Proposition~\ref{CorExpander}.
So, we have to prove that $D$
is an expander diagramm. 

Choose  $i\in \{0,1,\dots\}$ such that,
with 
$$
D_1= F\cup \bigcup_{j=1}^s\{ x_0^j,\dots, x_i^j\},
$$
we have $\mu(D_1) > 1/2.$

%$$D_0= F\bigcup \{x_0^1, x_0^2,\dots, x_0^s \}.$$
%For this, it suffices to show that,
%for every ray $r_j=\{x_0,x_1,\dots, \}$,
%there a positive lower for ${\mu (E(S, S^c))}/{\mu(S)}$
%uniformly for  sets $S$ contained in $\{x_1,x_2,\dots, \}$.

Let $S$ be  a subset of $D$  with $\mu(S)\leq \mu(D)/2.$
Then $D_1\nsubseteq S.$
Two cases can occur.

\noindent
$\bullet$\emph{First case:} $S\cap D_1=\emptyset.$ Thus, 
 $S$ is  contained in 
$$\bigcup_{j=1}^s\{x_{i+1}^j,x_{i+2}^j,\dots \}.$$
Fix $j\in\{1,\dots, s\}.$   
%and set for simplicity $r_j= \{x_0, ,x_1, x_2, \dots \}.$ 
%set $S_j=S\cap \{x_1^j,x_2^j,\dots, \}$
Let $i(j) \in \{0,1,\dots\}$ be minimal with the
property that $x_{i(j)+1}^j\in S.$ Then $e_j:=(x_{i(j)+1}^j,x_{i(j)}^j)\in E(S, S^c).$
Observe that $|\Ga_{x_{l+1}^j}|= {\rm deg}( x_l^j)|\Ga_{x_{l}^j}|$ for all $l\geq 0.$
Let $k$ be the minimal degree for vertices in $X$ (so, $k=\min\{k_0, k_1\}$
if $X= X_{k_0, k_1}$).Then $\mu(x_{l+1}^j) \leq\mu(x_{l}^j)/k$ for all $l$
and 
$$\mu(e_j)= \dfrac{1}{|\Ga_{e_j}|}\geq \dfrac{k}{|\Ga_{x_{i(j)}^j}|}=k\mu(x_{i(j)}^j).
$$
%Set $c=\mu (x_{i_0})=\dfrac{1}{|\Ga_{x_{i_0}}|}.$
Therefore, we have 
\begin{align*}
 \frac{\mu (E(S, S^c))}{\mu(S)}&\geq  \frac{\sum_{j=1}^s\mu(e_j)}{\sum_{j=1}^s\mu(\{x_{i(j)+1}^j,x_{i(j)+1}^j,\dots, \})}\\
%&\geq \sum_{j=1}^s\frac{\mu(e_j)}{\mu(\{x_{i(j)},x_{i+1},\dots, \})}\\
&\geq k \frac{\sum_{j=1}^s\mu(x_{i(j)}^j)} {\sum_{j=1}^s\sum_{l=0}^\infty \mu(x_{i(j)+l}^j)}\\
&\geq k \frac{\sum_{j=1}^s\mu(x_{i(j)}^j)}{\sum_{j=1}^s\mu(x_{i(j)}^j)\sum_{l=0}^\infty k^{-l}}\\
%&=\sum_{j=1}^s\frac{k}{{1}/{(1-k^{-1}}) }= s({k-1}).\\
&=k\frac{\sum_{j=1}^s\mu(x_{i(j)}^j)}{\frac{1}{1-k^{-1}}\sum_{j=1}^s\mu(x_{i(j)}^j)}\\
&=k\frac{1}{\frac{1}{1-k^{-1}}}=k-1.
\end{align*}

\noindent
$\bullet$\emph{Second case:} $S\cap D_1\neq\emptyset.$ 
Then there exist $x\in S\cap D_1$ and $y\in D_1\setminus S.$
Since $D_1$ is a connected subgraph, there exists a path
$(e_1,e_2,\dots, e_n)$ in $ED_1$ from $x$ to $y.$
Let $l\in\{1,\dots,n\}$ be minimal with the property $\partial_0(e_l)\in S$
and $\partial_1(e_l)\notin S.$ Then $e_l\in E(S,S^c)$.
Hence, with  $C=\min \{\mu(e)\ : \ e\in ED_1)\} >0,$  we have
$$\frac{\mu (E(S, S^c))}{\mu(S)}\geq \frac{C}{\mu(D)}.$$
This completes the proof of Theorem~\ref{Theo1}.$\bsq$

\section{Proof of Theorem~\ref{Theo2}}
\label{S4}
Let $(D,i,\mu)$ be a $k$-regular diagram. 
By the ``inverse    Bass--Serre theory'' of groups acting on trees,
there exists  a lattice $\Ga$ in  $G={\rm Aut} (X_k)$
for which $D=\Ga\backslash X_k$. Indeed, 
we can find a finite grouping of $(D,i)$, that is,
a graph of finite groups ${\bf D}=(D, \cal D)$ such that 
$i(e)$ is the index of ${\cal D}_e$ in ${\cal D}_{\partial_0 e}$
for all $e\in ED.$  Fix an origin $x_0.$ 
Let $\Ga=\pi_1({\bf D},x_0)$  be the fundamental group
of  $({\bf D},x_0)$. The universal covering of 
 $({\bf D},x_0)$ is the $k$-regular tree $X_k$ and 
the diagram $D$ can identified
with the diagram associated to  $\Ga\backslash X_k.$
For all this, see (2.5), (2.6) and (4.13) in \cite{BL}.

In view of Proposition~\ref{CorExpander},  Theorem~\ref{Theo2} will be proved once we  present
examples of $k$-regular diagrams with finite volume which are not expanders.
An example of such a diagram appears in \cite[Example 3.4]{Mokhtari}.
For the convenience of the reader, we review the construction.

Fix $k\geq 3$ and let $q=k-1.$
For every integer $n\geq 1$, let $D_n$ be the finite
graph with $2n+1$ vertices:
$$
\underset{x_1^{(n)}}\circ -\underset{x_2^{(n)}}\circ-\circ - \cdots \circ-\underset{x_{2n}^{(n)}}\circ-\underset {x_{2n+1}^{(n)}}\circ
$$
Let  $D$ be the following  infinite ray:
$$
\underset{x_0}\circ -\underset{x_1}\circ-D_1-\underset{x_2}\circ-\underset{x_3}\circ-D_2-\circ-\circ-\cdots- -\underset{x_{2n-2}}\circ-\underset{x_{2n-1}}\circ-D_n-\circ-\circ\cdots
$$
We first define a weight function $i_n$ on $ED_n$ as follows:
\begin{itemize}
 \item $i_n(e)=1$ if $e=(x_1^{(n)},x_2^{(n)})$ or  $e=(x_2^{(n)},x_1^{(n)})$
\item $i_n(e)= q$ if $e=(x_m^{(n)},x_{m+1}^{(n)})$ for $m$ even
\item $i_n(e)= 1$ if $e=(x_m^{(n)},x_{m+1}^{(n)})$ for $m$ odd
\item $i_n(e)= q$ if $e=(x_{m+1}^{(n)},x_m^{(n)} )$ for $m$ even
\item $i_n(e)= 1$ if $e=(x_{m+1}^{(n)},x_m^{(n)} )$ for $m$ odd.
\end{itemize}
Observe that $i_n(e)/i_n(\overline e)=1$ for all $e\in ED_n.$
Define now a weight function $i$ on $ED$ as follows:
\begin{itemize}
 \item $i(e)=q+1$ if $e=(x_0,x_1)$ 
\item $i(e)=q$ if $e=(x_1,x_0)$ 
\item $i(e)= 1$ if $e=(x_m,x_{m+1})$ for $m\geq 1$
\item $i(e)= q$ if $e=(x_{m+1},x_m)$ for $m\geq 1$
\item $i(e)= i_n(e)$ if $e\in ED_n$.
\end{itemize}

One readily checks that, for every vertex $x\in D,$ 
$$
\sum_{e\in \partial_0^{-1}(x)} i(e)= q+1 =k,
$$
that is,  $(D,i)$ is $k$-regular.
The  measure $\mu: VD\to \mathbf R^+$ corresponding to $i$
(see the remark at the beginning of Section~\ref{S2}) is given by
\begin{itemize}
 \item $\mu(x_0)=1/(q+1)$ 
\item $\mu(x_{2m-2})=1/q^{m-1}$ for  $m\geq 2$
\item $\mu(x_{2m-1})=1/q^m$ for  $m\geq 1$
\item $\mu(x)=1/q^n $ if $x\in D_n$.
\end{itemize}
One checks that, if we define
$\mu(e)= i(e) \mu(\partial_0 e)$ for all $e\in ED,$
we have $\mu(\overline e)= \mu(e).$ Moreover,
$$
\mu (D_n)=(2n+1)\frac{1}{q^n}
$$
and hence
$$
\mu(D)\leq \frac{1}{q+1} + 2\sum_{n \geq 0 }\frac{1}{q^n}+\sum_{n \geq 1 }\mu (D_n) <\infty.
$$
We have also
$$
E( D_n, D_n^c)= \{(x_{2n-1}, x_{2n-2}),(x_{2n}, x_{2n+1}) \}, 
$$
so that
$$
\mu\left(E( D_n, D_n^c)\right)= q\frac{1}{q^n} + \frac{1}{q^n}=\frac{q+1}{q^n}.
$$
Hence
$$
\frac{\mu\left(E( D_n, D_n^c)\right)}{\mu(D_n)}=\frac{\frac{q+1}{q^n}}{(2n+1)\frac{1}{q^n}}=
\frac{q+1}{2n+1}
$$
and 
$$\lim_n \frac{\mu\left(E( D_n, D_n^c)\right)}{\mu (D_n)} =0.$$
Observe that, since  $\lim_n\mu (D_n)=0,$ we have $\mu(D_n)\leq \mu(D)/2$
for sufficiently large $n.$
This completes the proof of Theorem~\ref{Theo2}.$\bsq.$

\noindent
{\bf Addresses}

\noindent
\obeylines
{Bachir Bekka
UFR Math\'ematiques, Universit\'e de  Rennes 1, 
Campus Beaulieu, F-35042  Rennes Cedex
 France}

\noindent
E-mail : bachir.bekka@univ-rennes1.fr

\medskip

\noindent
\obeylines
{Alexander Lubotzky
Institute of Mathematics, 
Hebrew University, Jerusalem 91904 
Israel}

\noindent
E-mail : alexlub@math.huji.ac.il

\end{document}